\documentclass{amsart}

\newtheorem{theorem}{Theorem}[section]
\newtheorem{lemma}[theorem]{Lemma}

\theoremstyle{definition}
\newtheorem{definition}[theorem]{Definition}
\newtheorem{example}[theorem]{Example}

\theoremstyle{corollary}
\newtheorem{corollary}[theorem]{Corollary}

\theoremstyle{remark}
\newtheorem{remark}[theorem]{Remark}

\numberwithin{equation}{section}

\def\Cal{\mathcal}

\def\H{{\Cal H}}

\def\R{{\Cal R}}

\def\F{{\Cal F}}

\def\I{{\Cal I}}

\def\Ma{\bbr^{n\times m}}

\def\Mt{\bbr^{k\times m}}
\def\Mkm{\bbr^{(n-k) \times m}}

\def\Mmm{\bbr^{m\times m}}

\def\Mnm{\bbr^{(n-m)\times m}}

\def \lv{{\bf \lam}}
\def \nv{{\bf n}_0}
\def \mn{{\bf m}_0}
\def \kv{{\bf k}_0}

\def\Z{\mathcal{Z}}

\def\f0{f_0}
\def\Fc0{\varphi_0}

\def\I_k {I_{-}^{k/2}}
\def\I+k {I_{+}^{k/2}}

\def\vnk{V_{n,k}}
\def\vnm{V_{n,m}}

\def\gnm{G_{n,m}}

\def\Gr{\frT}

\def\bbr{{\Bbb R}}

\def\bbc{{\Bbb C}}

\def\rank{{\hbox{\rm rank}}}

\def\const{{\hbox{\rm const}}}

\def\det{{\hbox{\rm det}}}

\def\min{{\hbox{\rm min}}}

\def\part{\partial}
\def\intl{\int\limits}

\def\Gam{\Gamma}
\def\Om{\Omega}
\def\a{\alpha}
\def\om{\omega}

\def\Del{\Delta}

\def\vp{\varphi}

\def\gam{\gamma}
\def\Lam{\Lambda}
\def\sig{\sigma}
\def\lam{\lambda}
\def\lv{{\boldsymbol \lam}}
\def\mv{{\boldsymbol \mu}}
\def\av{{\boldsymbol \a}}

\def\t{\tau}
\def\eq{\xi 'x=t}

\font\frak=eufm10

\def\fr#1{\hbox{\frak #1}}

\def\frL{\fr{L}}

\def\frT{\fr{T}}

\def\const{{\hbox{\rm const}}}

\def\det{{\hbox{\rm det}}}

\def\min{{\hbox{\rm min}}}

\def\p{\Om}
\def\gm{\Gamma_m}

\def\part{\partial}
\def\intl{\int\limits}

\def\Gam{\Gamma}
\def\Om{\Omega}
\def\a{\alpha}

\newcommand{\be}{\begin{equation}}
\newcommand{\ee}{\end{equation}}
\newcommand{\bea}{\begin{eqnarray}}
\newcommand{\eea}{\end{eqnarray}}
\newcommand{\Bea}{\begin{eqnarray*}}
\newcommand{\Eea}{\end{eqnarray*}}
\begin{document}

\title[ Composite Cosine Transforms] { Composite Cosine Transforms}

\author{ E. Ournycheva}
 \address{ Department of Mathematical Sciences, Kent State
 University,
 Mathematics and Computer
Science Building,
 Summit Street,
 Kent OH 44242,
  USA} \email{ournyche@math.kent.edu}

\thanks{The first author was supported by Abraham and Sarah Gelbart Research Institute for
Mathematical Sciences.}

\author{ B. Rubin}
 \address{Department of Mathematics, Louisiana State
University,  Baton Rouge, LA, 70803, USA, \hbox{\rm and} \newline
 Institute of Mathematics, Hebrew University, Jerusalem 91904,
 ISRAEL}
\email{borisr@math.lsu.edu}
\thanks{Both  authors were  supported in part
by the Edmund Landau Center for Research in Mathematical Analysis
and Related Areas, sponsored by the Minerva Foundation (Germany).}


\dedicatory{Dedicated to Professor
 Rolf Schneider on the occasion of
his 65th birthday}


\begin{abstract}
The cosine transforms of functions on the unit sphere play an
important role
   in convex geometry, the Banach space theory, stochastic geometry and other
   areas. Their higher-rank generalization to
 Grassmann manifolds represents an interesting mathematical object useful for applications.
   We introduce  more general  integral transforms
   that reveal distinctive features of higher
 rank objects in full generality. We call these new transforms  the
 {\it  composite cosine transforms}, by
taking into account that their kernels agree with the composite
power function of the cone of positive definite symmetric matrices.
 We show that injectivity of the composite cosine transforms can be
 studied using standard tools of the Fourier analysis on matrix spaces. In the framework of this
 approach, we introduce associated  generalized zeta
 integrals and give new simple proofs to the relevant functional
 relations. Our  technique   is based on  application of the higher-rank Radon transform on
matrix spaces.
 \end{abstract}

\maketitle

\section{Introduction}
\setcounter{equation}{0}

Let $S^{n-1}$ be the unit sphere in $\bbr^n$, $u \cdot v$ the usual
inner product of vectors $u, v \in S^{n-1}$ . The classical cosine
transform  \be\label{t11}(T f)(u)=\intl_{S^{n-1}} f(v) |u \cdot v|
\, dv, \qquad u \in S^{n-1},\ee and its generalization
\be\label{t1}(T^{\lam} f)(u)=\intl_{S^{n-1}} f(v) |u \cdot v|^{\lam}
\, dv,\ee are commonly in use in convex geometry, the Banach space
theory,  harmonic analysis, and many other areas; see
 \cite{Ga}, \cite{GH1}, \cite {Ko},   \cite{Schn}. Basic properties
 of $T^{\lam}$ (injectivity, boundedness in function spaces, and others)
 can be derived from the Funk-Hecke
 formula \be\label{sh}
T^{\lam}P_k=c\,\mu_k (\lam)\, P_k,\ee \be\label{mu1}
c=2\pi^{(n-1)/2}(-1)^{k/2},\qquad \mu_k (\lam)=
 \frac{\Gam\Big(\frac{\lam+1}{2}\Big)\,
\Gam\Big(\frac{k-\lam}{2}\Big)}{\Gam\Big(-\frac{\lam}{2}\Big) \,
\Gam\Big(\frac{\lam+k+n}{2}\Big)}, \ee where $P_k (x)$ is the
restriction of a homogeneous harmonic polynomial  of even degree $k$
\cite{Ru2}.

In the last two decades a considerable attention  was attracted to
higher-rank generalizations of  $T$ and $T^{\lam}$ for functions on
the Grassmann manifold $G_{n,m}$  of $m$-dimensional linear
subspaces of $\bbr^n$. We recall that, if  $\eta\in G_{n,m}$, $ \xi
\in G_{n,l}, \; l \ge m$, and $[\eta |\xi ]$ is  the $m$-dimensional
volume of the parallelepiped spanned by the orthogonal projection of
a generic orthonormal  frame in $\eta$ onto $\xi$, then, by
definition, \be\label{ccgr} (T^{\lam} f)(\xi)=\intl_{G_{n,m}}
f(\eta) \, [\eta | \xi]^{\lam} \, d\eta\ee (we adopt the same
notation $T^{\lam}$ as in (\ref{t1})). For $l>m$, the operator
(\ref{ccgr}) represents the composition of the similar one over
$G_{n,m}$ and the corresponding Radon transform acting from
$G_{n,m}$ to
 $G_{n,l}$ \cite {A}, \cite{GR}. Thus injectivity of $T^\lam$  can be
 studied using known results for the Radon transform
 (see \cite{GR} and references therein) and the case $l=m$ in (\ref{ccgr}).
 Owing to this, in the following we
assume that $l=m$, because just this case bears the basic features
of the operator family (\ref{ccgr}).

The investigation of operators (\ref{ccgr}) for $\lam=1$ was
initiated in stochastic geometry  by Matheron \cite[p. 189]{Mat1},
(see also \cite{Mat2}),
 who conjectured  that the higher-rank cosine transform $T^1$ is injective
as well as its rank-one prototype (\ref{t11}). Matheron's conjecture
was disproved in the remarkable paper by Goodey and Howard
\cite{GH1}. The higher-rank cosine transforms arise in convex
geometry in the context of the generalized Shephard problem for
lower dimensional projections \cite{GZ}. More general operators
$T^\lam$ for $\lam = 0,1,2,\dots$ were studied  in \cite[p.
117]{GH2}, where, by using reduction to $G_{4,2}$, it was proved
that $T^\lam$ is non-injective for such $\lam$. The range of the
$\lam$-cosine transform was characterized by Alesker and Bernstein
\cite{AB} ($\lam=1$) and by Alesker \cite{A} (any complex $\lam$),
in terms of representations of the special orthogonal group $SO(n)$.

In this  article we suggest a new  approach to  operators
$T^{\lam}$. This can be regarded as a complement to the well-known group
representation method.  The latter has proved
 to be  useful in the study of  Radon and cosine transforms,
 invariant differential and integro-differential operators on diverse
 homogeneous spaces of the orthogonal group; see \cite{Goo},
 \cite{Gr},  \cite{GH},  \cite{Str1}, \cite{Str2},
 \cite{TT}, and references therein. Our method differs
 from those in the cited papers. It gives
a direct analog of the multiplier equality (\ref{sh}) and is
applicable to a much more general operator family of the so-called
{\it composite cosine transforms}. These are introduced
 in Section 3  which describes main results of the paper. In Section
 4 we introduce the so-called {\it generalized zeta
integrals} with additional ``angle component" $f$   on the relevant
Stiefel manifold.   An important by-product of our investigation
 is  a
 functional equation  for these integrals that gives rise to the
 composite cosine  transform
  $T^{\lv}f$, for $\lv \in \bbc^m$,   in the most general form.
    The case $f\equiv 1$ was studied in \cite{FK} in the context
  of Jordan algebras. The argument from \cite{FK} was extended  in
   \cite{Cl} when
  $f$ is a determinantally homogeneous
 harmonic polynomial.  An alternative
proof of this result  for zeta integrals on matrix spaces was given
in \cite{OR}. This proof employs an idea from \cite{Kh2} to derive
the result for $\lam \in \bbc$ as a diagonal case  of the more
general statement for vector-valued $\lv =(\lam_1, \ldots \lam_m)
\in \bbc^m$. This idea allows us to avoid essential technical
difficulties (e.g., implementation of Bessel functions of matrix
argument) which arise when we get stuck on  the complex analysis of
a single variable; cf. \cite{FK}, \cite{Cl}, \cite{Ru5}.

In the present paper, we use the same idea and suggest  a new method
that demonstrates application of the higher-rank Radon transform on
matrix spaces \cite{OR1}, \cite{OR2}.
 This Radon transform enables
us to reduce the problem to the known case $f\equiv 1$. For the
usual cosine transform on the unit sphere, this approach  is due to
A. Koldobsky \cite{Ko}. Section 5 contains proofs of the main
results.

One should note that the  Fourier analysis of homogeneous
distributions is one of the oldest topics in the theory of
distributions, and there is a vast literature on this subject; see,
e.g., \cite{Es},\cite{GS}, \cite {Le}, \cite{Ra}, \cite{Sa},
\cite{Se}.

{\bf Acknowledgement.} We are thankful to Professors S. Alesker,
P. Goodey and W. Weil
 for very helpful comments and pleasant discussion.

\section{Preliminaries}

\setcounter{equation}{0}

The main references for this section are \cite{FK}, \cite{Gi},
 \cite{T}.
\subsection{Notation}
  Let $\bbr^{n\times m}$ be the
space of real matrices $x=(x_{i,j})$ having $n$ rows and $m$
 columns; $dx=\prod^{n}_{i=1}\prod^{m}_{j=1}
 dx_{i,j}$.
 In the following,
  $x'$ denotes the transpose of  $x$, $I_m$ is the identity $m \times m$
  matrix, and  $O(n)$   is the group of
real orthogonal $n\times n$ matrices. For $n\geq m$, we denote by
$\vnm= \{v \in \bbr^{n\times m}: v'v=I_m \}$   the Stiefel manifold
of orthonormal $m$-frames in $\bbr^n$. This is a homogeneous space
with respect to the action $\vnm \ni v\to \gam v $, $\gam\in O(n)$,
so that   $\vnm=O(n)/O(n-m)$.  The
  invariant measure $dv$ on $\vnm$ induced by the Lebesgue measure on the ambient space
is defined up to a constant multiple. We normalize it using
geometric argument and set
   \be \sigma_{n,m}
 \equiv \intl_{\vnm} dv =\prod\limits_{i=1}^{m} |S^{n-i}|, \ee
 where $|S^{i}|=2\pi^{(i+1)/2}/\Gam((i+1)/2)$ is the surface area of the $i$-dimensional unit
 sphere.

 Let  $\p$ be the cone
of positive definite symmetric matrices $r=(r_{i,j})_{m\times m}$
with the elementary volume $ dr=\prod_{i \le j} dr_{i,j}$. We denote
$|r|=\det(r)$ and let $ d_{*} r = |r|^{-(m+1)/2} dr$ be the
$GL(m,\bbr)$-invariant
 measure on $\p$. If $T_m$ is the group  of upper triangular $m \times m$
  matrices $t=(t_{i,j})$ with positive
diagonal elements, then each $r \in \p$ has a   unique
representation $r=t't$.

We will  constantly use the polar coordinates and the spherical
coordinates on $\bbr^{n\times m}$  which are defined as follows.
\begin{lemma}\label{l2.3} {\rm (\cite[pp. 66, 591]{Mu}, \cite{Ma})} If $x \in \Ma, \; \rank (x)=m, \; n \ge m$,
then \be \label{pol} x=vr^{1/2}, \qquad v \in \vnm,   \qquad r=x'x
\in\p,\ee and $dx=2^{-m} |r|^{(n-m-1)/2} dr dv$.
\end{lemma}
\begin{lemma}\label{sph} {\rm (\cite{P}, \cite{Ru5})} If
$x \in \Ma, \; \rank (x)=m,\;  n \ge m$, then
  \[
x=ut, \qquad u \in \vnm,   \qquad t \in T_m,\] and
$$
dx=\prod\limits_{j=1}^m t_{j,j}^{n-j} \,dt_{j,j}\,dt_*dv, \qquad
dt_*=\prod\limits_{i<j} dt_{i,j}.
$$
\end{lemma}

The Schwartz space  $S=S(\bbr^{n\times m})$ is
 identified with  the respective space on $\bbr^{nm}$ of infinitely differentiable rapidly decreasing functions.
 The Fourier transform  of a
function $f\in L^1(\bbr^{n\times m})$ is defined by \be\label{ft}
(\F f)(y)=\intl_{\bbr^{n\times m}} e^{{\rm tr(iy'x)}} f (x)
dx,\qquad y\in\bbr^{n\times m} \; .\ee The relevant Parseval
equality reads \be\label{pars} (\F f, \F \vp)=(2\pi)^{nm} \, (f,
\vp), \qquad (f, \vp)=\intl_{\bbr^{n\times m}} f(x)
\overline{\vp(x)} \, dx.\ee

\subsection{The composite power function}
Given $r=(r_{i,j})\in\p$, let $\Del_0(r)=1$, $\Del_1(r)=r_{1,1}$,
$\Del_2(r)$, $\ldots$, $\Del_m(r)=|r|$ be the corresponding
principal minors which are strictly positive. For
$\lv=(\lam_1,\dots,\lam_m)\in\bbc^m$, the composite power function
of the cone  $\p$  is defined by
 \be\label{pf} r^{\lv}=\prod\limits_{i=1}^m
 \left[\frac{\Del_i (r)}{\Del_{i-1} (r)}\right]^{\lam_i/2}=
  \Del_1 (r)^{\frac{\lam_1 -\lam_2}{2}}  \ldots \Del_{m-1}
(r)^{\frac{\lam_{m-1} -\lam_m}{2}} \Del_m (r)^{\frac{\lam_m}{2}}.
\ee
   In the special case $
\lam_1=\ldots=\lam_m=\lam$ we write $ \lv_0= (\lam, \ldots, \lam) \;
(\in \bbc^m) $ so that $r^{{\lv_0}}=|r|^{\lam/2}$. If $r=t't, \; t
=(t_{i,j})\in T_m$, then $r^{{\lv}}=\prod_{j=1}^m
 t_{j,j}^{\lam_j}$. This implies the following equalities:
\bea\label{pr1} r^{{\lv+\mv}}&=&r^{{\lv}}\;r^{{\mv}}, \quad
r^{{\lv+\av_0}}=r^{\lv}|r|^{\a/2}, \quad \av_0=(\a,\dots,\a);
\\ \label{pr6} \; (t'rt)^{{\lv}}&=&(t't)^{\lv}\;r^{{\lv}},\quad t\in
T_m. \eea The reverses of $\lv=(\lam_1,\dots,\lam_m)$ and
$r=(r_{i,j})\in\p$ are defined by \[ \lv_\ast=(\lam_m,\dots,\lam_1);
\qquad r_\ast =\om r\om,\qquad  \om=\left[\begin{array}{ccccc}
0 & {}   & {}    & 1 \\
                              {} & {}  & {.}    & {} \\
                              {} & {.}   & {}   & {} \\
                               1 & {}   & {}    & 0
\end{array} \right],\]
so that $$ (\lv_\ast)_j=\lam_{m-j+1}, \qquad
(r_\ast)_{i,j}=r_{m-i+1, m-j+1}.$$ We have
 \be\label{pr4}
 r^{{\lv_\ast}}=(r^{-1})_\ast^{-\lv},\qquad
(r^{-1})^{{\lv}}=r_\ast^{{-\lv_\ast}}.\ee

The  gamma function of the cone $\p$ is defined by \be\label{gf}
\Gam_{\Omega}  (\lv) =\intl_{\Omega} r^{\lv} e^{-{\rm tr} (r)} d_{*}
r=\pi^{m(m-1)/4}\prod\limits_{j=1}^{m} \Gam ((\lam_j- j+1)/2);\ee
 see, e.g., \cite[p. 123]{FK}. The integral in (\ref{gf}) converges
 absolutely if and only if $Re \, \lam_j>j-1$ for all $j=1,\dots,m$, and extends
 meromorphically to all $\lam \in \bbc^m$.
 The following relation holds:
 \be\label{eq11} \intl_{\Omega}
r^{\lv} e^{-{\rm tr} (rs)} d_{*} r= \Gam_{\Omega} (\lv) \,
s_\ast^{{-\lv_\ast}}, \qquad s\in\p.\ee
An important particular case
of (\ref{gf}) is  the Siegel integral \be\label{sgf} \Gam_{m} (\lam)
=\intl_{\Omega} |r|^{\lam} e^{-{\rm tr} (r)} d_{*}
r=\pi^{m(m-1)/4}\prod\limits_{j=0}^{m-1} \Gam (\lam - j/2), \quad Re
\,\lam
> (m-1)/2.\ee If $\lv_0=(\lam,\dots,\lam)$, then
 $ \Gam_{\Omega}  (\lv_0)=\gm(\lam/2)$. The volume $\sigma_{n,m}$ of the Stiefel manifold
 $\vnm$ may be written in terms of the Siegel Gamma function:   \be\label{2.16} \sigma_{n,m}
 = \frac {2^m \pi^{nm/2}} {\gm
 (n/2)}. \ee

 \subsection{Radon transforms on the space of matrices}
The main references for this subsection are  \cite{OR1}, \cite{OR2},
\cite{OR4}, \cite{P}, \cite{Sh1}, \cite{Sh2}. We fix  positive
integers $k,n$, and $m$, $0<k<n$, and let $\vnk$ be the Stiefel
manifold of orthonormal $k$-frames in $\bbr^n$. For $\; \xi\in\vnk$
and $t\in\Mt$, the linear manifold \be\label{plane} \tau=
\tau(\xi,t)=\{x\in\Ma:\eq\} \ee
 is called a {\it  matrix $(n-k)$-plane} in $\Ma$.  We denote by  $\Gr$ the
 set of all such
 planes.  Each $\t \in \Gr$ is
an ordinary $(n-k)m$-dimensional plane in $\bbr^{nm}$, but the set
$\frT$ has measure zero in the manifold of all such planes. The {\it
matrix Radon transform} $f(x)\to (\R_k f)(\t)$ assigns to a function
$f(x)$ on $\Ma$ a collection of integrals of $f$ over all matrix
planes $ \tau \in \Gr$,
 namely, \[ (\R_k
f) (\tau)=\int_{x \in \tau} f(x).\]  Precise meaning of this
integral is
 the following:
\be\label{4.9} (\R_k f) (\tau) =\intl_{\Mkm} f\left(g_\xi
\left[\begin{array} {c} \om \\t
\end{array} \right]\right)d\om,
\ee where $ g_\xi \in SO(n)$ is a rotation satisfying
\be\label{4.24}
g_\xi\xi_0=\xi, \qquad \xi_0=\left[\begin{array} {c}  0 \\
I_{k} \end{array} \right] \in \vnk. \ee

 The following statement is a matrix generalization of
the so-called projection-slice theorem. It links together the
Fourier transform (\ref{ft}) and the Radon transform (\ref{4.9}).
In the case $m=1$, this theorem can be found in [Na, p. 11] (for
$k=1$) and [Ke, p. 283] (for any $0<k<n$).
 \begin{theorem}\label{CST} {\rm(\cite{OR4})}
 For  $f\in L^1(\Ma)$ and  $1\le m\le k$,
\be\label{4.20} (\F f)(\xi b)=[\tilde \F (\R_k f)(\xi,\cdot)](b),
\quad \xi\in\vnk, \quad b\in\Mt. \ee Here, $\tilde \F\vp$ denotes
the Fourier transform of a function $t \to \vp(\xi, t)$ on the space
$\Mt$.
 \end{theorem}

\section{Main results}

 For $\lv=(\lam_1,\ldots, \lam_m)\in\bbc^m$, we introduce  intertwining operators \be
\label{tnf}(T^{\lv} f)(u)=\intl_{\vnm} f(v) \, (u'v v'u)^{\lv} \,
dv, \qquad u\in\vnm,\qquad n>m, \ee that commute with the left
action of $ \, O(n)$. We call $T^{\lv} f$ {\it the
 composite cosine transform  of } $f$. If
 $\lam_1=\ldots=\lam_m=\lam$, then
 (\ref{tnf}) reads
 \be\label{tll}(T^{\lam}
f)(u)=\intl_{\vnm} f(v) |\det( v' u)|^{\lam} \, dv. \ee If $f $ is
a $O(m)$ right-invariant function on
 $\vnm$, then (\ref{tll}) can be identified with (\ref{ccgr}) (for $l=m$) and
represents the usual $\lam$-cosine transform on $\gnm$.
\begin{definition}\label{def1} We denote by $\frL$  the set of all $\lv=(\lam_1, \dots ,\lam_m) \in
\bbc^m$ satisfying $Re \, \lam_j >j-m-1$ for all $j=1, \dots , m$.
\end{definition}
This definition is motivated by the following.
\begin{theorem}\label {exi}  For $f \in L^1(\vnm)$, the integral $(T^{\lv} f)(u)$
converges absolutely a.e. on  $\vnm$ if and only if $\lv \in
\frL$, and represents an analytic function of  $\lv$ in this
domain. For such $\lv$, the linear operator $T^{\lv} $ is bounded
on $L^1(\vnm)$. \end{theorem}

This statement follows immediately by Fubini's theorem from the
equality \cite{OR}
 \be \label {ave} \intl_{\vnm} (u'v v'u)^{\lv} \,
du=\frac{2^m \, \pi^{nm/2}}{\Gam_m (m/2)}\, \frac{\Gam_\Om(\lv
+\mn)}{\Gam_\Om(\lv+\nv)}\quad (\equiv T^{\lv} 1) \ee which is of
independent interest.

It is  challenging  to describe  the set of all $\lv \in \bbc^m$ for
which $T^{\lv}$ is injective. We cannot solve this problem in full
generality and restrict our consideration to the space $L^\flat
(\vnm)$ of $O(m)$ right-invariant integrable functions on $\vnm$.
This allows us to obtain   a precise description of those $\lv$ for
which $T^{\lv}$ is injective in the following important cases (a)
$2m \le n, \; $$ \lv=(\lam_1,\ldots, \lam_m)\in \bbc^m$, and (b)
$\lam_1=\dots =\lam_m=\lam \in \bbc$, provided that $T^{\lv}f$ and
$T^{\lam}f$ exist  as absolutely convergent integrals. Note that for
$m=1$, the space $L^\flat (\vnm)$ is actually the space of even
integrable functions on the unit sphere.
\begin{theorem}\label {main} Let $n> m\ge 1$ and $\lv \in \frL$, i.e.
 $Re \, \lam_j >j-m-1$ for all $j=1, \dots , m$.
 If, moreover, \be \label {uhh}\lam_j +m-j \neq 0, 2, 4,\dots
\quad \text{\rm for all} \quad j=1, \dots , m,\ee then the composite
cosine transform $T^{\lv}$ is injective on $L^\flat (\vnm)$. If
$2m\le n$ and (\ref{uhh}) fails, then $T^{\lv}$ is non-injective.
Specifically, it annihilates all $O(m)$ right-invariant, harmonic,
determinantally homogeneous polynomials of degree
$k>\max\limits_j \{\lam_j +m-j\}$ (see
Definition \ref{def2}).
\end{theorem}
Some comments are in order. The essence of our approach is that we
  apply the standard Fourier transform technique  to obtain
 a higher-rank analog of (\ref{sh}). We do not know if the condition
(\ref{uhh}) is necessary for injectivity of $T^\lv$ when $2m>n$. To
answer this question, one has to treat $T^\lv$ on polynomial
representations of $SO(n)$. These are parameterized by highest
weights $(m_1, m_2, \dots , m_{[n/2]})$ that are
more general than those adopted
 in the framework of our approach; cf. \cite{Str1}, \cite{TT}. However, if $\lam_1=\dots = \lam_m=\lam$, then
for the $\lam$-cosine transform (\ref{ccgr}), we  give the following
complete answer which reveals  essential difference between the
rank-one case and that of a higher-rank.
\begin{theorem}\label {main2} Let $n> m$, $Re\, \lam >-1$, and let
$r_{n,m}=\rank (\gnm)=\min (m, n-m)$ be the rank of the
Grassmannian $\gnm$. If $r_{n,m}=1$, then  $T^\lam$
 is injective on
$L^1 (\gnm) \;$ if and only if $\lam \neq 0,2,4, \dots \,$. If
$r_{n,m}>1$, then  $T^\lam$
 is injective on
$L^1 (\gnm) \;$ if and only if $\lam \neq 0,1,2, \dots \,$.
\end{theorem}

This statement   is known. It follows from the more general result
of Alesker \cite{A}.

\section{The generalized zeta integrals  and the composite cosine transforms}

\subsection{The generalized zeta integrals}
 In accordance with the polar decomposition $x=vr^{1/2},\; v \in
\vnm,   \; r=x'x \in\p$,   we introduce {\it the generalized zeta
integrals} (or zeta distributions): \bea\label{zeta}
\Z(\phi,\lv,f)&=&\intl_{\bbr^{n\times m}}
r^{\lv} \, f(v) \, \overline{\phi(x)} \,  dx=(r^\lv f, \phi), \\
\label{zeta*}\Z_*(\phi,\lv,f)&=&\intl_{\bbr^{n\times m}} r_*^{\lv}
\, f(v) \, \overline{\phi(x)}  \, dx=(r_*^\lv f, \phi), \eea where
$f \in L^1(\vnm)$ and $\phi \in S(\bbr^{n\times m})$. Zeta integrals
of this type with the angle component $f\equiv 1$ are well known and
arise in different occurrences; see \cite{FK}, \cite{BSZ}) and
references therein. We denote \bea\label{abs} {\bf
\Lam}&=&\{\lv\in\bbc^m :
Re\,\lam_j>j-n-1 \quad \text{\it for all} \quad j=1,\dots, m \}, \\
\label{pset}
{\bf \Lam}_0&=&\{\lv\in\bbc^m : \lam_j=j-n-l \quad \text{\it for some}\\
 &{}&j \in \{1,\dots, m\}, \quad \text{\rm and } \quad l \in \{1,3,5,\dots \} \}.
 \nonumber
 \eea
\begin{lemma}\label{l3.1}  The integrals (\ref{zeta}) and
(\ref{zeta*}) are  absolutely convergent if and only if $\lv \in
{\bf\Lam}$, and extend
 as  meromorphic functions of $ \lv$
with the  polar set ${\bf \Lam}_0 $.  The normalized zeta integrals
\be\label{nzi} \Z^0(\phi,\lv,f)=\frac{
\Z(\phi,\lv,f)}{\Gam_\Om(\lv+\nv)}, \qquad
\Z^0_*(\phi,\lv,f)=\frac{\Z_*(\phi,\lv,f)}{\Gam_\Om(\lv+ \nv)}, \ee
$\nv=(n, \ldots, n)$, are   entire functions of $\lv$.
 \end{lemma}

To prove this lemma it suffices to write both integrals in spherical
coordinates and then apply a standard argument  from \cite{GSh}; see
\cite{OR} for details.

\subsection{The basic functional equation}

The connection between zeta integrals and composite cosine
 transforms can be established  in the form of a
  functional equation which is actually the usual  Parseval
 equality.  Note that the function $(T^{\lv} f)(u)$, initially
 defined for $u \in \vnm$, extends to all matrices $y \in \bbr^{n\times m}$
of rank $m$. Indeed, by (\ref{pr6}), in spherical coordinates
$y=ut$, $ u \in \vnm$,   $ t \in T_m$, we have \be
\label{tf1}(T^{\lv} f)(y)=r^{\lv}(T^{\lv} f)(u)\ee where
$r^{\lv}=(t't)^{\lv}=(y'y)^{\lv}$ is the ``radial part" of $(T^{\lv}
f)(y)$.
\begin{theorem}\label{gen}
Let \be\label{fil} \vp_{\lv}(x)= \frac{ r_*^{-\lv_*-\nv}}{
\Gam_\Om(-\lv_*)} \, f(v), \qquad x=vr^{1/2}, \quad \nv=(n,\ldots,
n), \quad \lv\in\bbc^m. \ee If $f $ is  an integrable $O(m)$
right-invariant function on
 $\vnm$, then \be\label{eq5}
(\F\vp_{\lv})(y)=\frac{c_{\lv}}{\Gam_\Om(\lv+\mn)} \, (T^{\lv}
f)(y), \qquad c_{\lv}=2^{-|\lv|}\pi^{m^2/2}/\sig_{m,m},\ee
 in the sense of $S'$-distributions. In other words, for each $\phi\in
 S(\bbr^{n\times m})$,
\be \label{eq9}\frac{c_{\lv}}{\Gam_\Om(\lv+\mn)} \, ( T^{\lv} f,\;
\F\phi)=(2\pi)^{nm} \, (\vp_\lv, \; \phi) \equiv (2\pi)^{nm} \,
\Z^0_*(\phi,-\lv_*-\nv, f).  \ee
\end{theorem}
A self-contained proof of this statement is given in \cite{OR}.
Here, we apply an alternative approach  which  is of independent
interest. The main idea is to reduce the problem to the
corresponding functional equation containing zeta integrals on
$\Mmm$ with the angle component $f\equiv 1$. To this end, we invoke
the higher-rank Radon transform (\ref{4.9}). We start with  two
auxiliary lemmas.
\begin{lemma}\label{l-kh} For $r\in\p$ and $y\in\Ma$,
\be \label{eq3}\F\Big[\frac{ r_*^{-\lv_*-\nv}}{
\Gam_\Om(-\lv_*)}\Big](y)=\frac{2^{-|\lv|}\pi^{nm/2}}{\Gam_\Om(\lv+\nv)}\;(y'y)^{\lv}
\ee in the sense  of $S'$-distributions.
\end{lemma}
\begin{proof} Formula (\ref{eq3}) was established by Khekalo
\cite{Kh2} who modified    the  argument from [St2, Chapter III,
Sec. 3.4] for functions of matrix argument; see also \cite{FK},
\cite{OR4}, \cite {Ru5} on this subject. For convenience of the
reader, we outline the proof of (\ref{eq3})  in our notation. Since
$$ \F[ e^{-{\rm tr}( xs x'/4\pi)}](y)=(2\pi)^{nm}|s|^{-n/2} e^{-{\rm
tr}(\pi
 y s^{-1}y')}, $$ for $s\in\p$ and $\phi\in S(\bbr^{n\times m})$,  the Parseval equality yields    \bea
\label{eq-21}&&|s|^{-n/2}\intl_{\Ma} e^{-{\rm tr}(\pi y s^{-1}y')}\,
\overline{(\F\phi)(y)}\, dy=\intl_{\Ma} e^{-{\rm tr}( xs x'/4\pi)}\,
\overline{\phi(x)}\, dx. \eea We multiply (\ref{eq-21})  by
$s^{\lv+\nv}$ and integrate against $d_\ast s$. This gives
\[ \intl_{\Ma} I_1(y)\, \overline{(\F\phi)(y)}\, dy=\intl_{\Ma} I_2(x)\, \overline{\phi(x)}\, dx, \]
 where
\[I_1(y)=\intl_{\p} s^{\lv} e^{-{\rm tr}(\pi y
s^{-1}y')}d_\ast s,\quad I_2(x)=\intl_{\p}s^{\lv+\nv } e^{-{\rm tr}(
xs x'/4\pi)}d_\ast s. \] Evaluation of  the last integrals by means
of (\ref{eq11}) gives
  $$ I_1(y)=\pi
^{|\lv|/2}\Gam_{\Omega} (-\lv_*) (y' y)^{ \lv}, \quad Re\,\lam_j
<j-m,
$$  and $$
I_2(x)=(4\pi)^{(|\lv|+nm)/2}\Gam_{\Omega} (\lv+\nv) (x' x)_*^{-\lv_*
-\nv},\quad Re\,\lam_j >j-n-1.$$ Hence, if $ j-n-1<Re\,\lam_j<j-m$,
then \bea \label{eq8} &&\Gam_{\Omega} (-\lv_*) \intl_{\Ma}  (y' y)^{
\lv} \,\overline{(\F\phi)(y)}\, dy=c_\lv\Gam_{\Omega} (\lv+\nv)
\intl_{\Ma}(x' x)_*^{-\lv_* -\nv} \,\overline{\phi(x)}\, dx, \eea
where $c_\lv=2^{nm+|\lv|}\pi^{nm/2}. $ By Lemma \ref{l3.1}, this
extends analytically to all $\lv$, and we are done.
\end{proof}
\begin{lemma}
For $\phi\in S(\Ma)$  and  $\xi\in\vnm$, \be\label{eq1}
\frac{1}{\Gam_\Om(\lv+\mn)}\intl_{\Ma}\!\!\! (x'\xi\xi ' x)^{\lv} \;
\overline{(\F \phi)(x)} dx= \frac{d_{\lv}}{\Gam_\Om(-\lv_*)}
 \intl_{\Mmm}\!\!\!(t't)_*^{-\lv-\mn}\overline{\phi(\xi t)}dt,\ee \be\label{clam2}
d_{\lv} =2^{mn+|\lv|}\pi^{nm-m^2/2}.\ee
\end{lemma}
\begin{proof}
Let us denote by $A(\xi)$ the left side of (\ref{eq1}). By passing
to the spherical  coordinates on $\Ma$, according to Lemma \ref{sph}
and (\ref{pr6}) we obtain
$$
A(\xi)=\frac{1}{\Gam_\Om(\lv+\mn)}\intl_{T_m} \prod\limits_{j=1}^m
t_{j,j}^{\lam_j+n-j} dt_{j,j}dt_*\intl_{\vnm}
\overline{(\F\phi)(ut)} \, (u'\xi \xi'u)^{\lv} \, du.
$$
By Theorem \ref {exi}, this integral converges absolutely if  $\lv
\in \frL$, that is  $Re \, \lam_j
>j-m-1$ for all $j=1, \dots , m$. Let us  evaluate $A(\xi)$
for such $\lv$.

We put $x=g_\xi \left[\begin{array} {c} \omega
\\t
\end{array} \right]$, where $ g_\xi \in SO(n)$ is a rotation satisfying
(\ref{4.24}) with $k=m$, $\om\in\Mnm$, and $t\in\Mmm$. Then
 $\xi'x=t$, so that the
Fubini theorem  and (\ref{eq3}) yield \bea\nonumber
A(\xi)&=&\frac{1}{\Gam_\Om(\lv+\mn)}
 \intl_{\Mmm} (t't)^{\lv} \; dt \intl_{\Mnm}\; \overline{(\F
\phi)(g_\xi \left[\begin{array} {c} \omega \\t
\end{array} \right])} d\om\\[14pt] \nonumber &=&\frac{1}{\Gam_\Om(\lv+\mn)}
 \intl_{\Mmm} (t't)^{\lv} \overline{(\R_m \F\phi)(\xi,
t)}\; dt \\[14pt] \nonumber
&=& \frac{2^{|\lv|+m^2} \pi^{m^2/2}}{\Gam_\Om(-\lv_*)}
 \intl_{\Mmm} (t't)_*^{-\lv_*-\mn} \overline{(\tilde\F^{-1}\R_m \F\phi)(\xi,
t)}\; dt\\[14pt] \nonumber
&=&\frac{2^{|\lv|} \pi^{-m^2/2}}{\Gam_\Om(-\lv_*)}
 \intl_{\Mmm} (t't)_*^{-\lv_*-\mn} \overline{(\tilde\F\R_m \F\phi)(\xi,
t)}\; dt. \eea Here,  $\tilde\F$ and $\tilde\F^{-1}$ denote the
Fourier transform acting in the $t$-variable on the space $\Mmm$,
and its inverse, respectively.  We have  used the equality
$(\tilde\F^{-1}\vp)(\xi, t)=(2\pi)^{-m^2}( \tilde\F\vp )(\xi, -t)$.
Moreover, we have applied (\ref{eq3}) with $n=m$ in the case when
both functionals in this equality are regular.  The latter is true
if  \be\label{r-lam} j-m-1<Re\,\lam_j<j-m. \ee By (\ref{4.20}),
$$
\tilde \F [\R_m \F\phi(\xi, \cdot)](t)=(\F \F\phi)(\xi
t)=(2\pi)^{nm}\phi(-\xi t).
$$
This gives \bea \nonumber A(\xi)&=& \frac{2^{|\lv|+nm}
\pi^{nm-m^2/2}}{\Gam_\Om(-\lv_*)}
 \intl_{\Mmm} (t't)_*^{-\lv_*-\mn} \overline{\phi (\xi
t)}\; dt. \eea Therefore, (\ref{eq1}) is valid when $\lv$ satisfies
(\ref{r-lam}). By Lemma \ref{l3.1}, this  extends analytically to
all $\lv\in\bbc^m$.
\end{proof}

{\bf Proof of Theorem \ref{gen}.} Let   $\lv$ obey (\ref{r-lam}), so
that by Theorem \ref{exi} and Lemma \ref{l3.1}, integrals  in
(\ref{eq9}) are absolutely convergent. Then, by (\ref{eq1}) and the
Fubini theorem,
 \bea \nonumber \frac{1}{\Gam_\Om(\lv+\mn)}( T^{\lv} f,\; \F\phi)&=&\frac{1}{\Gam_\Om(\lv+\mn)}\intl_{\Ma}(T^{\lv} f)(x)
 \overline{(\F  \phi)(x)} dx\\
 \nonumber&=&\frac{1}{\Gam_\Om(\lv+\mn)}
 \intl_{\vnm} f(v) dv\intl_{\Ma}\overline{(\F  \phi)(x)}  (x'v v' x)^{\lv}  \,dx\\\nonumber &=&\frac{d_{\lv}}{\Gam_\Om(-\lv_*)}\intl_{\vnm} f(v) dv
 \intl_{\Mmm}(t't)_*^{-\lv_*-\mn}\overline{\phi(v t)}dt.\eea
Here, $d_{\lv}$ is defined by (\ref{clam2}).  By the polar
decomposition (see Lemma \ref{l2.3}),  for the right-invariant
function $f$  we obtain \bea \nonumber \frac{1}{\Gam_\Om(\lv+\mn)}(
T^{\lv} f,\; \F \phi)\!\!\!&=&\!\!\!
\frac{2^{-m}d_{\lv}\sig_{m,m}}{\Gam_\Om(-\lv_*)}\!\!\!\intl_{\vnm}\!\!\!
f(v) dv\!\intl_{O(m)}\!\!\!d\gam
 \intl_{\p}r_*^{-\lv_*-\mn}|r|^{-1/2}\overline{\phi(v \gam r^{1/2})}dr  \\\nonumber
&=& \frac{2^{-m}d_{\lv}\sig_{m,m}}{\Gam_\Om(-\lv_*)}\intl_{\vnm}
f(v) dv
 \intl_{\p}r_*^{-\lv_*-\mn}|r|^{-1/2}\overline{\phi(v  r^{1/2})}dr
 \nonumber \\
&=& \frac{d_{\lv}\sig_{m,m}}{\Gam_\Om(-\lv_*)}\intl_{\Ma}
f(x(x'x)^{-1/2}) (x'x)_*^{-\lv_*-\nv}\overline{\phi(x)}dx.\nonumber
 \eea
Hence,  (\ref{eq9}) is proved  for $j-m-1<Re\,\lam_j<j-m$. According
to Lemma \ref{l3.1}, this extends analytically to all $\lv\in
\bbc^m$.
\begin{example} Let $m\ge 1, \; \lam_1=\ldots=\lam_m=\lam, \; |x|_m
=\det(x'x)^{1/2}$. Then
\[\vp_{\lv}(x)= \frac{|x|_m^{-\lam-n}}{ \Gam_m (-\lam/2)} \,
f(x(x'x)^{-1/2}), \qquad (T^{\lam} f)(y)=\intl_{\vnm} f(v)
|\det(v' y)|^{\lam} \, dv. \] If $f $ is  an integrable $O(m)$
right-invariant function on
 $\vnm$, then the operator $T^{\lam}$ coincides with (\ref{ccgr}),   and we have \be\label{rnm}
(\F\vp_{\lv})(y)=\frac{2^{-\lam m} \, \pi^{m^2/2}
}{\sig_{m,m}\Gam_m((\lam +m)/2)} \, (T^{\lam} f)(y).\ee If $m=1$,
and  $x, y \in \bbr^n \setminus \{0\}$, then
\[\vp_{\lv}(x)= \frac{|x|^{-\lam-n}}{ \Gam (-\lam/2)} \,
f\Big(\frac{x}{|x|}\Big), \qquad (T^{\lam} f)(y)=\intl_{S^{n-1}}
f(v) |v \cdot y|^{\lam} \, dv. \]  In this case, \be\label{rn1}
(\F\vp_{\lv})(y)=\frac{2^{-1-\lam}\pi^{1/2}}{\Gam((\lam +1)/2 )} \,
(T^{\lam} f)(y).\ee The last equality is well known and can be found
in many sources.
\end{example}

\subsection{The case of  homogeneous
polynomials}

Let $P_k(x)$ be a polynomial  on $\bbr^{n\times m}$ which is
harmonic (as a function on $\bbr^{nm}$) and determinantally
homogeneous of degree $k$, i.e., $ P_k(xg)=\det (g)^k P_k(x),\quad
\forall g\in GL(m,\bbr)$. It means that $P_k$ is a usual homogeneous
harmonic polynomial of degree $km$ on $\bbr^{nm}$. Theorem \ref{gen}
can be
 strengthened  if we choose $f$ to be the restriction of
$P_k (x)$ onto $\vnm$.
\begin{theorem} Let
\be\label{fik} \vp_{\lv, k}(x)= \frac{ r_*^{-\lv_*-\nv}}{
\Gam_\Om(-\lv_*)} \, P_k(v), \qquad x=vr^{1/2}, \quad
d_{\lv}=2^{-|\lv|}\pi^{nm/2}i^{km}.\ee Then, for all $\lv \in
\bbc^m$, \be\label{eq50} (\F\vp_{\lv,k})(vr^{1/2})=\frac{d_{\lv}\,
\Gam_\Om(\kv-\lv_*)}{\Gam_\Om(-\lv_*) \, \Gam_\Om(\lv+\kv+\nv)} \,
P_k(v)\, r^{\lv}\ee in the sense of $S'$-distributions. In other
words, for each $\phi\in S(\bbr^{n\times m})$, \bea \frac{d_{\lv}\,
\Gam_\Om(\kv-\lv_*)}{\Gam_\Om(\lv+\kv+\nv)} \,( P_k(v)\, r^{\lv},\;
\F\phi)&=&(2\pi)^{nm}(\vp_{\lv,k}, \; \phi) \nonumber
\\&=& (2\pi)^{nm}\,\Z^0_*(\phi,-\lv_*-\nv, P_k).\nonumber \eea
\end{theorem}
\begin{proof} The classical Hecke identity
$$ \int_{\bbr^{n\times m}}P_k(x) e^{-{\rm tr}(\pi x'  x)}e^{{\rm tr}(2\pi i
y'x)}\;dx=i^{km}P_k(y)e^{-{\rm tr}(\pi y'  y)}$$ implies that \bea
\label{eq-2}&& \!\!\!\!|s|^{-k-n/2}\intl_{\bbr^{n\times m}}
P_k(y)e^{-{\rm tr}(\pi y s^{-1}y')}\, \overline{(\F\phi)(y)}\,
dy\\\nonumber &=&(2\pi i)^{-km}\intl_{\bbr^{n\times m}} P_k(x)
e^{-{\rm tr}( xs x'/4\pi)}\, \overline{\phi(x)}\, dx. \eea
 We
multiply  this by $s^{\lv+\nv+\kv}$ and  proceed, as in the proof of
Lemma \ref{l-kh}, to obtain the result for $
j-n-k-1<Re\,\lam_j<j+k-m$. Since this domain is not void for all
$k=0,1, 2, \ldots$, and the normalized zeta integral
$\Z^0_*(\phi,-\lv_*-\nv, P_k)$ is  an entire function of $\lv$, the
result follows by analytic continuation.\end{proof}
\begin{example} Let $m\ge 1, \; \lam_1=\ldots=\lam_m=\lam$, $\; |x|_m
=\det(x'x)^{1/2}$,
$$\vp_{\lam, k}(x)=\frac{|x|_m^{-\lam-n}}{ \Gam_m (-\lam/2)} \, P_k
(x(x'x)^{-1/2}), \qquad d_{\lam}=2^{-\lam m}\pi^{nm/2}i^{km}.$$ Then
\be \label{ft2} (\F\vp_{\lam,k})(y)=\frac{ d_{\lam} \,
\Gam_m((k\!-\!\lam)/2)}{\Gam_m(-\lam/2) \,
\Gam_m((\lam\!+\!k\!+\!n)/2)} \, |y|_m^{\lam} \, P_k
(y(y'y)^{-1/2}).\ee If $m=1$, then $\vp_{\lam,
k}(x)=\frac{|x|^{-\lam-n}}{ \Gam (-\lam/2)} \, P_k
\Big(\frac{x}{|x|}\Big)$, and we have \be \label{ft1}
(\F\vp_{\lam,k})(y)\!=\!\frac{d_{\lam} \,
\Gam((k-\lam)/2)}{\Gam(-\lam/2) \, \Gam((\lam+k+n)/2)} |y|^{\lam}
P_k \Big(\frac{y}{|y|}\Big),\quad
d_{\lam}\!=\!2^{-\lam}\pi^{n/2}i^{k}.\ee
\end{example}

  \section{Proofs of  main results}

 \subsection{A higher-rank analog of (\ref{sh})}
\begin {definition}\label{def2} \cite{Herz} A polynomial $P_k(x)$
on $\bbr^{n\times m}$ is called {\it an $H$-polynomial of degree
$k$} if it is $O(m)$ right-invariant, harmonic, and determinantally
homogeneous of degree $k$. We denote by $\H_k$ the space of all such
polynomials. \end{definition}
\begin{lemma} \label{inff} Let $P_k \in \H_k$,
\be\label{mul} \mu_k(\lv)=\frac{\Gam_\Om(\lv+\mn)\,
\Gam_\Om(\kv-\lv_*)}{\Gam_\Om(-\lv_*) \, \Gam_\Om(\lv+\kv+\nv)},
\qquad \lv\in \bbc^m.\ee If $\lv$
 does not belong to the polar set of $\Gam_\Om(\lv+\mn)$, then
\be\label{cmp} (T^{\lv} P_k)(vr^{1/2})=c \,\mu_k(\lv)\,
P_k(v)r^{\lv}, \qquad c=\pi^{m(n-m)/2} \, i^{km}\,\sig_{m,m},\ee in
the sense of $S'$-distributions.
\end{lemma}
\begin{proof} Let us  compare (\ref{eq5}) and (\ref{eq50}), assuming  $f(v)=P_k(v)$.
 For all $\lv \in
\bbc^m$,  \be\frac{c_{\lv}}{\Gam_\Om(\lv+\mn)} \, (T^{\lv}
P_k)(vr^{1/2})
 =\frac{d_{\lv}\,
\Gam_\Om(\kv-\lv_*)}{\Gam_\Om(-\lv_*) \, \Gam_\Om(\lv+\kv+\nv)} \,
P_k(v)\, r^{\lv},\ee
\[ c_{\lv}=2^{-|\lv|}\pi^{m^2/2}/\sig_{m,m}, \qquad
d_{\lv}=2^{-|\lv|}\pi^{nm/2}i^{km},\] in the sense of $
S'$-distributions. If we exclude all $\lv$ belonging to the polar
set of $\Gam_\Om(\lv+\mn)$,  we get (\ref{cmp}).
\end{proof}
\begin{corollary} \label{cr} $\text{\rm (Cf. (\ref{sh}))}$ \  If $\lv \in \frL$
(see Definition \ref{def1}) and $P_k \in \H_k$, then
\be\label{lmu}(T^{\lv} P_k)(v)=c \,\mu_k (\lv)\, P_k(v), \qquad v
\in \vnm,\ee $c$ and $\mu_k (\lv)$ being the same  as in (\ref{cmp}).
\end{corollary}
\begin{proof}  The function $\Gam_\Om(\lv+\mn)$ has no poles in
$ \frL$. Hence,  by (\ref{cmp}), we have $ ((T^{\lv} P_k)(vr^{1/2}),
\phi)=c \,\mu_k(\lv)\, (P_k(v)r^{\lv}, \phi) $ for all $\phi
(y)\equiv \phi (vr^{1/2}) \in S(\bbr^{n\times m})$. Choose $\phi
(y)=\chi (r)\, \psi (v)$, where $\chi (r)$ is a non-negative
$C^\infty$ cut-off function supported away from the boundary of $\p$
and $\psi (v)$ is a $C^\infty$
 function on $\vnm$. By passing to polar coordinates, owing to
   (\ref{tf1}) we obtain
 $$c_\chi \,\intl_{\vnm} [(T^{\lv}
P_k)(v)-c \,\mu_k(\lv)\,P_k (v)]\, \psi (v)\, dv =0, \qquad
c_\chi=\const \neq 0.$$ This implies (\ref{lmu}).
\end{proof}
\begin{remark} An important question is, do there exist $H$-polynomials of a given
 degree $k$? For $n=m$ we have exactly two such
polynomials, namely, $P_0(x)\equiv 1$ and $P_1(x)=\det (x)$. It is
known \cite[p. 484]{Herz} that for $2m\le n$ there exist
$H$-polynomials of every degree $k$.  The space $\H_k$ in this case
is spanned by polynomials of the form $P_k(x)=\det(a'x)^k$ where $a$
is a complex $n\times m$ matrix satisfying $a'a=0$ \cite[p. 27]{TT}.
\end{remark}

\subsection{Proof of Theorem \ref{main} }
 To prove the first
statement, we consider the equality \be \label {fr1}((T^{\lv}
f)(vr^{1/2}),\; \F\phi)=A(\lv) \; (r_*^{-\lv_*-\nv} \, f(v), \;
\phi),\ee $$A(\lv)=(2\pi)^{nm}\,
\Gam_\Om(\lv+\mn)/c_{\lv}\,\Gam_\Om(-\lv_*),$$ which follows from
(\ref{eq9}) and (\ref{fil}). Suppose that  $(T^{\lv} f)(v)=0$ a.e.
on $\vnm$ for some $\lv \in \frL$. Then $(T^{\lv} f)(y),\;
y=vr^{1/2} \in \bbr^{n\times m}$, is zero for almost all $y \in
\bbr^{n\times m}$, and (\ref{fr1}) yields $A(\lv)\;
(r_*^{-\lv_*-\nv} \, f(v), \; \phi)=0$. The assumption (\ref{uhh})
along with $\lv \in \frL$ imply $Re \, \lam_j
>j-m-1 $ and $ \lam_j \neq j-m, j-m+2,\dots \, $.
Hence $\lv$ is not a pole of $\Gam_\Om(-\lv_*)$, and therefore
$A(\lv) \neq 0$. This gives $(r_*^{-\lv_*-\nv} \, f(v), \; \phi)=0,$
where  the left  side is understood in the sense of analytic
 continuation. Choosing $\phi$ as in the proof of Corollary
 \ref{cr}, we obtain $f(v) =0$ a.e. on $\vnm$.

To prove the second statement, we  note that, for $2m<n$,
$H$-polynomials of every degree $k$ do exist.  We observe that the
function $\Gam_\Om(\kv-\lv_*)$ in  $\mu_k(\lv)$ can be written as \[
\Gam_\Om(\kv-\lv_*)=\pi^{m(m-1)/4}\prod\limits_{j=1}^{m} \Gam
\Big(\frac{k+j-\lam_j- m}{2}\Big). \] It has no poles in $\frL$ if $
k>\max\limits_j \{\lam_j +m-j\}$.
 Since $\Gam_\Om(\lv+\mn)$ also has no poles in $\frL$, then, by
(\ref{lmu}), $T^\lv P_k =0$ for all such $k$ provided that $\lv$ is
a pole of $\Gam_\Om(-\lv_*)$ (i.e.,  (\ref{uhh}) fails). This proves
the theorem.

\subsection{Proof of Theorem \ref{main2} }
 We make use of the canonical homeomorphism of $G_{n,m}$ and
$G_{n,n-m}$ under which $T^\lam $ is invariant. To be precise, given
$u, v \in \vnm$, let $\xi\in\gnm$ and $\eta\in \gnm$ be
$m$-dimensional subspaces of $\bbr^n$ spanned $u$ and $v$,
respectively. We denote by $\xi^\perp\in G_{n,n-m}$ and
$\eta^\perp\in G_{n,n-m}$ the corresponding orthogonal subspaces,
and choose  any $ u_\perp \in V_{n, n-m}$ in $\xi^\perp$ and $
v_\perp \in V_{n, n-m}$ in $\eta^\perp$. Since $[\eta | \xi]=[\eta
^\perp | \xi^\perp]$, we can successively define the functions
$F(\eta)$ on $\gnm$, $F_\perp(\eta^\perp)$ on $G_{n,n-m}$, and
$f_\perp(v_\perp)$ on $V_{n, n-m}$ by $ F(\eta)=f(u), \;
F_\perp(\eta^\perp)=F((\eta^\perp)^\perp), \;
f_\perp(v_\perp)=F_\perp(\eta^\perp)$. Then $T^\lam $ expresses
through the similar operator $T_{\perp}^\lam$ on  $V_{n, n-m}$ as
follows: \be\label{expr} (T^\lam f)(u)=(T_{\perp}^\lam
f_\perp)(u_\perp).\ee Indeed,
 \bea(T^\lam f)(u)&=&\intl_{G_{n,m}} F(\eta) \, [\eta |
\xi]^{\lam} \, d\eta=\intl_{G_{n,n-m}}F_\perp(\eta^\perp)\,[\eta
^\perp | \xi^\perp]^\lam \, d\eta^\perp\nonumber \\
&=&\intl_{V_{n,n-m}} f_\perp(v_\perp) |\det(v_{\perp}'
u_\perp)|^{\lam} \, dv_\perp=(T_{\perp}^\lam
f_\perp)(u_\perp).\nonumber\eea If $\rank (G_{n,m})=1$, i.e., $m=1$
or $n-1$, the result follows from the multiplier equality (\ref{sh})
which is a particular case of (\ref{lmu}).

Let $\rank (G_{n,m})>1$. Owing to (\ref{expr}), it suffices to prove
the theorem for $2m\le n$, because otherwise we have $2(n-m)<n$, and
one can treat  $T_{\perp}^\lam$ instead of $ T^\lam$.   If
$\lam_1=\ldots =\lam_m=\lam$, the condition $\lam \neq 0,1,2, \dots$
coincides with (\ref{uhh}). But  for $2m\le n$,  the result follows
from
 Theorem \ref{main}.

\bibliographystyle{amsalpha}

\begin{thebibliography}{A}

\bibitem [A] {A} S. Alesker, \textit{The $\a$-cosine transform and
intertwining integrals}, Preprint, 2003.

\bibitem [AB] {AB} S. Alesker, and J. Bernstein, \textit{Range
characterization of the cosine transform on higher Grassmannians},
Advances in Math., \textbf{184} (2004), 367--379.

\bibitem [BSZ] {BSZ} L. Barchini, M. Sepanski, and R. Zierau, \textit{Positivity of Zeta
distributions and small representations}, Preprint, 2004.

\bibitem [Cl] {Cl}
 J.-L. Clerc, \textit{Zeta distributions associated to a
representation of a Jordan algebra}, Math. Z. \textbf{239} (2002),
 263--276.

\bibitem [Es] {Es}  G. I. Eskin,  \textit{ Boundary value problems for elliptic
pseudodifferential equations}, Amer. Math. Soc., Providence, R.I.,
1981.

\bibitem [FK] {FK} J. Faraut, and A. Kor\'anyi, \textit{Analysis on
symmetric cones}, Clarendon Press, Oxford, (1994).

\bibitem [Ga] {Ga} R.J. Gardner, \textit{ Geometric tomography}, Cambridge University
Press, New York, 1995.



\bibitem [G\v{S}] {GS}  Gel'fand, I.M., and  \v{S}apiro, Z.Ja., \textit{ Homogeneous
functions and their applications}, Uspekhi Mat. Nauk, \textbf{10},
(1955), no. 3, 3--70 (in Russian).

\bibitem [GSh] {GSh} I. M. Gel'fand, and G. E. Shilov,  \textit{Generalized functions}, Vol. 1. \textit{ Properties
and operations},
 Academic Press, New
York-London, 1964.

\bibitem [Gi]{Gi} S.G. Gindikin, \textit{Analysis on homogeneous
domains}, Russian Math. Surveys, \textbf{19}  (1964), No. 4,
1--89.

\bibitem [GH] {GH} F. Gonzalez, and S. Helgason,  \textit{Invariant differential
operators on Grassmann manifolds}, Adv. in Math. \textbf{60}
(1986),  81--91.

\bibitem [Goo] {Goo} P. Goodey, \textit{Applications of representation
 theory to convex bodies}, II International Conference in "Stochastic Geometry,
 Convex Bodies and Empirical Measures" (Agrigento, 1996),
   Rend. Circ. Mat. Palermo (2) Suppl.  No. 50 (1997), 179--187.

\bibitem  [GH1] {GH1} P. Goodey, and R. Howard, \textit{Processes
of flats induced by higher-dimensional processes}, Adv. in Math.,
\textbf{80} (1) (1990), 92--109.

\bibitem  [GH2] {GH2} P. Goodey, and R. Howard, \textit{Processes of
flats induced by higher-dimensional processes. II.} Integral
 geometry and tomography (Arcata, CA, 1989), 111--119, Contemp. Math., 113,
 Amer. Math. Soc., Providence, RI, 1990.

\bibitem  [GZ] {GZ} P. Goodey, and G. Zhang, \textit{Inequalities between
projection functions of convex bodies},  Amer. J. Math. \textbf{120}
(1998), 345--367.



\bibitem [Gr]{Gr}  E. L. Grinberg,  \textit{Radon transforms on higher
Grassmannians}, J. Differential Geom. \textbf{24} (1986), 53--68.

\bibitem [GR] {GR}  E. Grinberg, and  B. Rubin,  \textit{Radon inversion on Grassmannians
 via G{\aa}rding-Gindikin fractional integrals}, Annals of Math.
 \textbf{159} (2004), 809--843.


\bibitem [Herz]{Herz} C. Herz, \textit{Bessel functions of matrix
argument}, Ann. of Math., \textbf{61} (1955), 474--523.

\bibitem [Kh] {Kh2}  S.P. Khekalo, \textit{The Igusa zeta function associated
with a composite power function on the space of rectangular
matrices}, Preprint POMI RAN, 10 (2004), 1--20.


\bibitem [Ko] {Ko} A. Koldobsky,   \textit{ Inverse formula for the
Blaschke-Levy representation}, Houston J. Math.,\textbf { 23}
(1997), 95-107.

\bibitem [Le] {Le} C. Lemoine,  Fourier transforms of homogeneous distributions,
 {\it Ann. Scuola Norm. Super. Pisa Sci. Fis. e Mat.}, {\bf 26}
(1972), No. 1, 117--149.

\bibitem [Ma] {Ma} A.M. Mathai, \textit{ Jacobians of matrix
transformations and functions of matrix argument}, World Sci.
Publ. Co. Pte. Ltd, Singapore, 1997.


\bibitem [Mat1] {Mat1} G. Matheron, \textit{Un th\'eor\`eme d'unicit\'e pour les
 hyperplans poissoniens}, J. Appl. Probability, \textbf {11}
 (1974), 184--189.

\bibitem [Mat2] {Mat2} \bysame, \textit{Random sets and integral geometry},
 Wiley Series in Probability and Mathematical Statistics. John Wiley
\& Sons, New York-London-Sydney, 1975.

\bibitem [Mu] {Mu} R.J. Muirhead, Aspects of multivariate
statistical theory, John Wiley \& Sons. Inc., New York, 1982.

\bibitem [OR1] {OR1}  E. Ournycheva,  and B. Rubin,  \textit{Radon transform of functions  of
matrix argument}, Preprint, 2004 (math.FA/0406573).

\bibitem [OR2] {OR}  E. Ournycheva,  and B. Rubin,  \textit{The Composite Cosine Transform on
 the Stiefel Manifold and Generalized Zeta Integrals}, Contemporary Math. (to appear).

\bibitem [OR3] {OR2}  E. Ournycheva,  and B. Rubin,  \textit{Higher-rank Radon
transforms}, Preprint 2005.

\vskip 1 truecm

\bibitem [OR4] {OR4}  E. Ournycheva,  and B. Rubin,  \textit{An analogue
of the Fuglede formula in integral geometry on matrix spaces},
Contemporary
 Math., \textbf {382} (2005), 305--320.


\bibitem [P] {P} E.E. Petrov,  \textit{The Radon transform in
spaces of matrices}, Trudy seminara po vektornomu i tenzornomu
analizu,  M.G.U., Moscow, \textbf{15} (1970), 279--315 (Russian).

\bibitem [Ra] {Ra} M. Ra\"{\i}s,  \textit{Distributions homog\`enes sur des
espaces de matrices}, Bull. Soc. math. France, Mem., \textbf{30},
(1972), 3--109.

\bibitem [Ru] {Ru2} B. Rubin, \textit{ Inversion of fractional integrals related
to the spherical Radon transform},  Journal of Functional Analysis,
\textbf{157} (1998), 470--487.

\bibitem [Ru4] {Ru5} \bysame, \textit{Riesz potentials and integral
geometry in the space of rectangular matrices},  Advances in Math.
 (to appear).

\bibitem [Sa] {Sa} Samko, S. G., \textit{ Generalized Riesz
potentials and hypersingular integrals with homogeneous
characteristics, their symbols and inversion},
  Proceeding of the Steklov Inst. of Math.,  \textbf{ 2} (1983) , 173--243.

\bibitem [Schn] {Schn}   R. Schneider, \textit{ Convex bodies: The Brunn-Minkowski theory},
Cambridge Univ. Press, 1993.

\bibitem [Se] {Se}  V.I. Semyanistyi,   \textit{Some integral transformations and integral
geometry in an elliptic space},   Trudy Sem. Vektor. Tenzor. Anal.,
\textit{12}  (1963), 397--441 (in Russian).

\bibitem [Sh1] {Sh1} L.P. Shibasov,   \textit{ Integral problems in a matrix space that are
connected with the functional $ X\sp{\lambda }\sb{n,m}$.}
 Izv. Vys\v s. U\v cebn. Zaved. Matematika (1973), No.
8 (135), 101--112 (Russian).


\bibitem [Sh2] {Sh2}  \bysame, \textit{Integral geometry on planes of a matrix space}.
(Russian) Harmonic analysis on groups. Moskov. Gos. Zao\v cn. Ped.
Inst. Sb. Nau\v cn. Trudov Vyp. \textbf{39}
   (1974), 68--76.



\bibitem [Str1] {Str1} R.S. Strichartz,  \textit{The explicit Fourier
decomposition of } $L^2(SO(n)/SO(n-m))$, Can. J. Math., \textbf{27}
(1975), 294--310.

\bibitem [Str2] {Str2} \bysame  \textit{$L\sp p$ estimates for Radon transforms in
Euclidean and non-Euclidean spaces}, Duke Math. J.  \textbf{48}
(1981),  699--727.

\bibitem [T] {T} A. Terras, \textit{Harmonic analysis on symmetric spaces
and applications}, Vol. II,  Springer, Berlin, 1988.

\bibitem [TT] {TT}  T. Ton-That, \textit{Lie group representations and harmonic
polynomials of a matrix variable},  Trans. Amer. Math. Soc. \textbf{
216} (1976), 1--46.

\end{thebibliography}

\end{document}